\documentclass[11pt,a4paper]{amsart} 
\usepackage{amscd, amsmath,amsthm,amsxtra}
\usepackage{amssymb}

\usepackage[utf8]{inputenc}
\usepackage[english]{babel}

\usepackage[all]{xy}

\newcommand{\cO}{\mathcal{O}}

\newcommand{\bC}{\mathbb{C}}

\newcommand{\N}{\mathbb{N}}

\newcommand{\Q}{\mathbb{Q}}
\newcommand{\R}{\mathbb{R}}

\newcommand{\Z}{\mathbb{Z}}

\newtheorem{Pro} {Proposition}
\newtheorem{lem}{Lemma}
\newtheorem{Cor}{Corollary}

\begin{document}
\title[LCK metrics on OT manifolds]{On Locally Conformally K\" ahler 
 metrics on Oeljeklaus-Toma Manifolds }
%\author{\c Stefan Deaconu, Victor Vuletescu}
\author{\c Stefan Deaconu }
\address{Faculty of Mathematics and Informatics\\
	University of Bucharest,
	Academiei st. 14,
	Bucharest, Romania}
\email{stefan.deaconu.al@gmail.com}

\author{Victor Vuletescu}
\address{Faculty of Mathematics and Informatics\\
	University of Bucharest,
	Academiei st. 14,
	Bucharest, Romania}
\email{vuli@fmi.unibuc.ro}



\thanks{Both authors are partially supported by Romanian Ministry of Education and Research,
Program PN-III, Project number PN-III-P4-ID-PCE-2020-0025, Contract 30/04.02.2021}
\date{} % no date}
\maketitle

\begin{abstract} 
We show that Oeljeklaus-Toma manifolds $X(K, U)$ where $K$ is a number field of signature $(s, t)$ such that $s\geq 1,$  $t\geq 2$  and $s\geq 2t$ admit no lck metric. Combined with the earlier results by \cite{OeTo} and \cite{Dub} this completely solves the problem of existence of LCK metrics on Oeljeklaus-Toma manifolds.
\end{abstract}
%\tableofcontents
\section{Introduction}
%\subsection{Oeljeklaus-Toma (OT) manifolds}
Oeljeklaus-Toma manifolds were introduced by  K.Oeljeklaus and M. Toma in \cite{OeTo}, as a generalization to higher dimensions of the Inoue surfaces $S_M$ (\cite{Inoue}). Very briefly, their construction goes as follows. Fix a number field $K$ having $s\geq 1$ real embeddings and $2t\geq 2$ complex ones, and label its embeddings such that $\sigma_1,\dots, \sigma_s$ are the real ones, while $\sigma_{s+t+i}=\overline{\sigma}_{s+i}$ for any $i=1,\dots, t.$ Let $\cO_K$ be the ring of integers of $K,$ $\cO_K^*$ the group of units of $\cO_K$ and $\cO_K^{*, +}$ the subgroup of $\cO_K^*$ of {\em totally positive} units, that is elements $u\in \cO_K^*$ such that $\sigma_i(u)>0$ for all $i=1,\dots, s.$ Letting ${\mathbb H}:=\{z\in \bC\vert Im(z)>0\}$, we see there are natural actions of $\cO_K$ and respectivley of $\cO_K^{*,+}$ on ${\mathbb H}^s\times \bC^t\subset \bC^{s+t}$ by
$$a\cdot (x_i)_{i=1,\dots s+t}:=\left(x_i+\sigma_i(a)\right)_{i=1,\dots s+t}, \; \forall a \in \cO_K$$
and respectively
$$u\cdot (x_i)_{i=1,\dots s+t}:=\left(\sigma_i(u)x_i\right)_{i=1,\dots s+t}, \; \forall u \in \cO_K^{*, +}.$$
The combined resulting action of $\cO_K^{*, +}\ltimes \cO_K$ is however not discrete in general. Still, in \cite{OeTo} it is shown that one can always find subgroups $U\subset \cO_K^{*, +}$ such that the action of $U\ltimes \cO_K$ is discrete and cocompact: the resulting compact complex manifold is usually denoted $X(K, U)$ and is called an {\em Olejeklaus-Toma manifold} (OT, for short). The Inoue surfaces $S_M$ are corresponding to the particular case when $K$ has degree $3$ and $U\subset \cO_K^{*, +}$ is any subgroup of finite index.

The OT manifolds have a number of very interesting properties. For instance,they are non-K\"alerian, do not respect Hodge symmetry (nor satisfy the $\partial\overline{\partial}$-lemma), and for appropriate (but rather generic) choices of the group of units $U$, they have no closed proper complex analytic subspaces.

%\subsection{Special metrics on non-K\"alerian manifolds}

As these manifolds do not admit K\" ahler metrics, it is natural to ask whether other natural metrics to exist on them: for a detailed account on this problem, see e.g.  \cite{ADD}. One of the most interesting ones are the {\em locally conformally K\"{a}hler} metrics (LCK, for short): these are those whose associated $(1,1)-$ forms $\omega$ has the property
\begin{eqnarray}\label{lck}
d\omega=\theta\wedge \omega
\end{eqnarray}
for some closed $1-$form $\theta$ (for more details see \cite{DO}).
The existence of such metrics on OT manifolds $X(K, U)$ can be read on the Galois properties of the group of units $U$. More precisely, it was shown (see \cite{Dub}, appendix by L. Battisti) that: 

\begin{Pro}\label{Bat}
An Oeljeklaus-Toma manifold $X(K, U)$ admits an LCK metric if  if and only if for any unit $u\in U$ one has
\begin{eqnarray}\label{rest}
\vert \sigma_{s+1}(u)\vert=\dots=\vert \sigma_{s+t}(u)\vert
\end{eqnarray} 
\end{Pro}

Already since these manifolds were introduced in \cite{OeTo}, it was shown that such metrics exists on $X(K, U)$ when $t=1$ and do not exist when $s=1$ and $t\geq 2.$
For the remaining possibilities for $(s,t)$, the second-named author showed in (\cite{Vu}) that for a certain number of cases, LCK metrics do not exist. The result was widely extended by A. Dubickas in \cite{Dub}; still, some cases remained unclear. The goal of this note is to show the non-existence of LCK metrics  in the remaining cases.

\medskip

{\small {\em Acknowledgments.} The authors thank Liviu Ornea for useful discussions and valuable suggestions  and Alexandra Otiman for a careful reading of a previous version of the present note.}

\section{The results}
\begin{lem}\label{lem}
Let $\Lambda \subset \R^n$ be a discrete lattice. Then $\Lambda$ can not be written as a union
\begin{eqnarray}\label{uni}
\Lambda=\bigcup_{i=1}^m\Lambda_i
\end{eqnarray}
of sublattices of proper rank, $rank_{\Z}(\Lambda_i)<rank_\Z(\Lambda), \forall i=1,\dots, m.$
\end{lem}

{\bf Proof.} For a lattice $\Lambda\subset \R^n$ freely generated by some vectors $e_1,\dots, e_N$ we let
$$\Lambda \Q:=\{\sum_{j=1}^Nq_ie_i, q_i\in \Q,\forall i\};$$ then $\Lambda \Q$ is a $\Q-$vector space and if $\Lambda$ is discrete then
$dim_{\Q}(\Lambda \Q)=rank_{\Z}(\Lambda).$

Next, we infer that (\ref{uni}) implies 
\begin{eqnarray}\label{uniq}
\Lambda\Q=\bigcup_{i=1}^m\Lambda_i\Q
\end{eqnarray}
The only inclusion to see here is $''\subset''.$ Take a vector $v\in \Lambda \Q;$ then 
$v=\displaystyle\sum_{j=1}^Nq_ie_i$ hence
$v=\frac{1}{M}\displaystyle\sum_{j=1}^Na_ie_i$ for some $M\in \N^*$ and $a_i\in \Z.$
But then the vector $w:= \displaystyle\sum_{j=1}^Na_ie_i$ is in $\Lambda$, hence by our assumption $w$ belongs to some $\Lambda_i$; it follows that $v\in \Lambda_i\Q$.

But decomposition (\ref{uniq}) leads to a contradiction, since a vector space over an infinite field cannot be written as a finite union of subspaces of smaller dimension. Q.E.D.

\medskip

{\em 
{\bf Notations.} Let $K$ be a number field with $s$ real embeddings and $2t$ complex ones; we suppose we labeled the embeddings $\sigma_i (i=1\dots s+2t)$ of $K$ such that $\sigma_1,\dots, \sigma_s$ are the real embeddings and such that $\overline{\sigma}_{s+k}=\sigma_{s+k+t}$ for all $k=1,\dots, t.$
We denote by $\Lambda_K$ the image of the units of $K$ under the logarithmic embedding
$$l(u):=\left(\log\vert \sigma_1(u)\vert,\dots , \log\vert \sigma_{s+t}(u)\vert\right)\subset \R^{s+t}.$$
Dirichlet units theorem tells us that $\Lambda_K$ is a discrete (and complete) lattice in the  hyperplane ${\mathcal H}_{Dir}$ given by
$$({\mathcal H}_{Dir}):\; x_1+\dots +x_s+2x_{s+1}+\dots +2x_{s+t}=0$$
If $L\subset K$ is a number subfield, we will similarly denote by $\Lambda_L$ the image of the units of $L$ under the previous embedding.%; more generally, if one has a subgroup of $\cO_K^*$ defined by a multiplicative relation, say $$({r}):\;\sigma_1(u)^{n_1}\dots \sigma_s(u)^{n_s}\sigma_{s+1}(u)^{n_{s+1}}\dots \sigma_{s+2t}(u)^{n_{s+2t}}(u)=1, \; (n_j\in \Z, for\; all\; j)$$we will denote by $\Lambda_r$ its image under the logarithmic  embedding, say$$(\Lambda_r):\; n_1x_1+\dots +n_sx_s+(n_{s+1}+n_{s+t+1})x_{s+1}+\dots +(n_{s+t}+n_{s+2t})x_{s+t}=0$$
}

\begin{Pro}\label{nolck} Let $K$ be a number field of signature $(s, t)$. If $s\geq 1,$  $t\geq 2$  and $s\geq 2t$ then $\cO_K^*$ has no subgroup $U$ of rank $s$ such that
$$\sigma_{s+1}(u)\sigma_{s+t+1}(u)=\dots=\sigma_{s+t}(u)\sigma_{s+2t}(u)$$
holds good for any $u\in U.$

\end{Pro}

{\bf Proof.} Assume such an $U$ would exist.
First notice that the logarithmic image $l(U)$ of $U$ lives on the intersection of the hyperplane ${\mathcal H}_{Dir}$ above with  the $t-1$ hyperplanes ${\mathcal H}_i, i=1, \dots, t-1$ given by 
$$x_{s+i}=x_{s+i+1}.$$
Notice that the linear variety ${\mathcal H}_{Dir}\cap\left(\displaystyle\bigcap_{i=1}^{t-1}{\mathcal H}_i\right)$ is of dimension $s.$
%The image $l(U)$ under the logarithmic embedding of $U$  is a lattice of rank $s.$ 

\medskip

We will prove that there are finitely many sublattices $\Lambda'\subset \Lambda_K$ with $rank(\Lambda')<s$ such that  any element $l(u)\in l(U)$ lives in (at least)  one such $\Lambda'$, getting henceforth a contradiction with the Lemma \ref{lem}. So take an arbitrary element $u\in U.$

If $deg(u)<[K:\Q]$ then there exists some proper subfield $L\subsetneq K$ such that $u\in L;$ in particular $l(u)\in \Lambda_L.$ Call $(s', t')$ the signature of $L$; then $rank(\Lambda_L)=s'+t'-1.$ Letting $d:=[K:L]$ we have $s+2t=d(s'+2t').$ Hence  $s'+2t'=\frac{s+2t}{d}$ so $$s'+t'-1=\frac{s+2t}{d}-t'-1.$$ Now $$\frac{s+2t}{d}-t'-1<s\Leftrightarrow s+2t<ds+d(t'+1)\Leftrightarrow 2t-d(t'+1)<(d-1)s.$$ But $2t-d(t'+1)<2t$ and $(d-1)s\geq s$ as $d\geq 2$ (since $L\subsetneq K$). Hence $rank(\Lambda_L)<s.$

We are left with the case when $u$ has maximal degree. 
As $t\geq 2$ the relation
 \begin{eqnarray}\label{rel}
 \sigma_{s+1}(u)\sigma_{s+t+1}(u)=\sigma_{s+2}(u)\sigma_{s+t+2}(u)
 \end{eqnarray}
 holds good. 
 %As $u$ has maximal degree and 
 As the absolute Galois group  $G_{\overline{\Q}/\Q}$  %of  the splitting field  of minimal polynomial of $u$ 
 acts transitively on the Galois conjugates of $\sigma_{s+1}(u)$, we see there exists some $\varphi\in  G_K$  such that $\varphi(\sigma_{s+1}(u))=\sigma_{1}(u).$
 Applying $\varphi$ to relation (\ref{rel}) we get 
 \begin{eqnarray}\label{rel2}
 \sigma_1(u)\sigma_j(u)=\sigma_k(u)\sigma_l(u)
 \end{eqnarray}
  for some $j, k, l\in \{2,\dots, s+2t\}.$
    Taking absolute values we see that the logarithmic  image of $u$ lives in the hyperplane ${\mathcal H}_{j'k'l'}$ of $\R^n$ given by 
    $$({\mathcal H}_{j'k'l'}):\; x_1+x_{j'}=x_{k'}+x_{l'}$$
    where $j', k', l'$ equals respectively $j, k, l$ if they are $\leq s+t$ or $j-t, k-t, l-t$ otherwise.
    Since $s\geq 1$ we see that
     $$dim({\mathcal H}_{j'k'l'}\cap {\mathcal H}_{Dir}\cap\left(\bigcap_{i=1}^{t-1}{\mathcal H}_i\right))=s-1<s$$
    so $$\Lambda_{j'k'l'}:=\Lambda_K\cap\left({\mathcal H}_{j'k'l'}\cap {\mathcal H}_{Dir}\cap\left(\bigcap_{i=1}^{t-1}{\mathcal H}_i\right)\right)$$
    is a lattice of rank $<s$ since it is discrete.
    
    We conclude  that any $l(u), u\in U$ lives   either in a lattice of the form $\Lambda_L$ with $L\subsetneq K$ a proper subfield or 
    in a lattice of the form 
    $\Lambda_{i'j'k'}$ as above; as all these lattices are of rank $<s$ and they are finitely many, we got our contradiction.
    
    \begin{Cor} If an Oeljeklaus-Toma manifold $X(K, U)$ admits an LCK metric, then $K$ has exactly $2t=2$ complex embeddings.
    \end{Cor}
    
    The proof follows at once from the above Propositions \ref{Bat} and \ref{nolck}.


\begin{thebibliography}{100}
\bibitem[AnDuOtSt22]{ADD} Angella, D., Dubickas, A., Otiman, A., Stelzig, J., {\em On metric an cohomological properties of Oeljeklaus-Toma manifolds}, ArXiv: 2201.06377
\bibitem[DrOr98]{DO} Dragomir,S.,   Ornea, L.,  {\em Locally conformally K\"{a}hler manifolds}, Progress in
Math. 55, Birkh\"{a}user, 1998.
\bibitem[Dub14]{Dub} Dubickas, A.,  {\em Nonreciprocal units in a number field with an application to Oeljeklaus–Toma manifolds}, New
York J. Math, 20(2014), p. 257–274.
\bibitem[Ino74]{Inoue} Inoue, M., {\em On surfaces of class V II0}, Invent. Math., 24(1974), 269–310
\bibitem[MilneANT]{carte} James S. Milne. {\em Algebraic Number Theory}. Available at
http://www.jmilne.org/math/CourseNotes/ant.html.

\bibitem[OeTo05]{OeTo} Oeljeklaus, K., Toma, M., {\em Non-K\"ahler compact complex manifolds associated to number
fields}, Ann. Inst. Fourier, Grenoble, 55(2005), no. 1, 161–171.
\bibitem[Vu14]{Vu}  Vuletescu V., {\em LCK metrics on Oeljeklaus-Toma manifolds versus Kronecker’s theorem}, Bull. Math. Soc. Sci.
Math. Roumanie, 57 (2014), no. 2, p. 225-231 
\end{thebibliography}
\end{document}